\documentclass[12pt]{article}
\usepackage{latexsym}
\usepackage{amssymb}
 
\newcommand{\proof}{{\noindent \bf Proof. }}
\newtheorem{thm}{Theorem}

\newtheorem{lem}{Lemma}
\newtheorem{cor}{Corollary}
\newtheorem{prop}{Proposition}
\newtheorem{conj}{Conjecture}
\newtheorem{rem}{Remark}

\newcommand\mbf[1]{\mbox{\boldmath$#1$}}
\newcommand\msbf[1]{\mbox{\boldmath\scriptsize$#1$}}
\newcommand{\N}{{\mathbb N}}

\date{}

\begin{document}
\begin{titlepage}
\title{\bf ON TYPES OF GROWTH FOR GRAPH--DIFFERENT PERMUTATIONS}
{\author{{\bf  J\'anos K\"orner}
\\{\tt korner@di.uniroma1.it}
\\''La Sapienza'' University of Rome
\\ ITALY
\and {\bf G\'abor Simonyi}\thanks{Research partially supported by the
Hungarian Foundation for Scientific Research Grant (OTKA) Nos.\ 
T046376, AT048826, NK62321, and EU ToK Project FIST 003006.}
\\{\tt simonyi@renyi.hu}
\\R\'enyi Institute of Mathematics, Budapest
\\ HUNGARY
\and {\bf Blerina Sinaimeri}
\\{\tt sinaimeri@di.unirma1.it}
\\''La Sapienza'' University of Rome
\\ ITALY}}

\maketitle
\begin{abstract}

We consider an infinite graph $G$ whose vertex set is the set of natural
numbers and adjacency depends solely on the difference between vertices. 
We study the largest cardinality of a set of permutations of $[n]$
any pair of which differ somewhere in a pair of adjacent vertices of $G$ and
determine it completely in an interesting special case. We give estimates for
other cases and compare the results in case of complementary graphs. We also
explore the close relationship between our problem and the concept of Shannon
capacity "within a given type". 
\end{abstract}
\end{titlepage}

\section{Introduction}

The topic of our paper has its origins in the following 
mathematical puzzle of K\"orner and Malvenuto \cite{KM}. Call two permutations
of $[n]:=\{1,\dots,n\}$ {\em colliding} if, represented by linear orderings of
$[n]$, they 
put two consecutive elements of $[n]$ somewhere in the same
position. For the maximum cardinality $\rho(n)$ of a set of 
pairwise colliding permutations of $[n]$ the following conjecture was
formulated.  

\begin{conj}\label{conj:clak} {\rm (\cite{KM})}
For every $n\in \N$ 
$$\rho(n)={n \choose {\left \lfloor \frac{n}{2} \right \rfloor}}.$$
\end{conj} 

This conjecture remains open; for the best bounds the interested reader may
consult  \cite{KMS} and \cite{BF}. 
In this paper we initiate a systematic study of similar problems for all the
graphs 
on the countable vertex set $\N$. Not only do we believe that these problems
are interesting on their own, but beyond this we hope that studying them  
within a unified framework may shed more light also on the initial problem on
colliding permutations. 

Let  $G$ be an infinite graph whose vertex set is the set $\N$ of the natural
numbers.  
We call two permutations of the elements of $[n]$, the first $n$ natural
numbers,  
$G$--different if they map some $i \in [n]$ to adjacent vertices of $G.$ 
(We will often think and write about permutations of $[n]$ as $n$-length
sequences that contain each element of $[n]$ exactly once. In this language,
two permutations are $G$--different if there exists a position where the 
corresponding two sequences contain the two endpoints of an edge of $G$.) 
We denote by 
$T_G(n)$ the maximum cardinality of a set of pairwise $G$--different
permutations  
of $G.$ The question about pairwise colliding permutations is that special
case of our  
present problem where the graph is the infinite path $L$ defined by
$$V(L)=\N \qquad E(L)=\{\{i, i+1\} \; | \; i\in \N\}.$$
We will concentrate our attention on the special class of {\it distance
  graphs}.  
Given an arbitrary (finite or infinite) set $D\subseteq \N$ we define the
graph  
$G=G(D)$ by setting
$$V(G(D))=\N \qquad E(G(D))=\{\{i, i+d\} \; | \; i\in \N, d\in D\}.$$
Clearly, $L=G(\{1\})$. 
We  will write 
$$T(n, D)=T_{G(D)}(n).$$

In the papers \cite{KM} and \cite{KMS} attention was restricted to those cases
where  
the growth of $T(n, D)$ in $n$ is only exponential. Here we consider various 
speeds of growth and determine 
the exact value of $T(n, D)$ for every $n$ in a non--trivial case.
We are especially interested in the relationship between the values of $T(n,
D)$ and $T(n, \overline{D})$, where $\overline{D}=\N\setminus D$. 

\section{Superexponential growth}

The determination of $T(n, \overline{\{1\}})$ leads to a surprisingly simple formula. 
\begin{thm}\label{thm:com}
$$T(n, \overline{\{1\}})={n! \over 2^{\lfloor{n \over 2}\rfloor}} \quad \hbox{for every} \quad n\in \N.$$

\end{thm}

\proof
First we prove the upper bound
$$T(n, \overline{\{1\}})\leq {n! \over 2^{\lfloor{n \over 2}\rfloor}}.$$
To this end fix $n$ and define $\sigma_{i,j}$ to be the permutation that
exchanges the entries $i\in [n]$ and $j\in [n]$, that is, for any permutation
$\pi$, $\sigma_{i,j}\pi$ differs from $\pi$ only in the places where the
entries $i$ and $j$ stand, which are exchanged.
For any fixed $\pi$ consider the set of permutations
$$C(\pi):=\{\sigma_{1,2}^{\varepsilon_{1,2}}\sigma_{3,4}^{\varepsilon_{3,4}}\dots\sigma_{k,k+1}^{\varepsilon_{k,k+1}}\pi:
\forall i\ \varepsilon_{i,i+1}\in\{0,1\}\},$$ 
where $k$ equals  $2\lfloor n/2\rfloor-1$,  
$\sigma_{i,j}^0$ is meant to be the identity permutation, while
$\sigma_{i,j}^1:=\sigma_{i,j}$. 
Let $B$ be a set of permutations of $[n]$ satisfying our condition that for
any pair of them there is an $i \in [n] $ they map to numbers at distance at
least two and observe that the conditions imply  $|C(\pi)\cap B|\leq 1$,
while $C(\pi)\cap C(\pi')=\emptyset$ if $\pi,\pi'\in B, \pi\neq\pi'$.  
Since $|C(\pi)|=2^{\lfloor n/2\rfloor}$ for any $\pi ,$ the foregoing implies
$$|B|\leq {n! \over 2^{\lfloor{n \over 2}\rfloor}},$$
which is the claimed upper bound.

In order to prove the inequality in the opposite direction, for every $n$ we
shall explicitely construct a set of permutations satisfying the requirement. 
We start by the odd values of $n$ and build our construction in a recursive manner. 
It will be important for the recursion that for every odd $n$ the construction be invariant with respect to cyclic shifts.
For $n=1$ the construction consists of the identical permutation. 
Suppose next to have constructed 
$$t_{n-2}:={(n-2)! \over 2^{\lfloor{(n-2) \over 2}\rfloor}}$$
permutations yielding a set $B_{n-2}$ that satisfies the pairwise relation we need 
and has the additional property of being closed with respect to cyclic shifts.
We will construct a set $A_n$ of ${n-1 \over 2}t_{n-2}$ permutations of $[n]$
satisfying the same pairwise condition and define $B_n$ to be the set
consisting of all the  
cyclically shifted versions of the elements of $A_n.$ For an arbitrary
permutation $\pi$ of $[n-2]$ and $1<j\leq n$ we define the transformations   
 $\Psi^j$ in the following manner. The permutation $\Psi^j\pi$ is acting on the set $[n]$,  
\begin{eqnarray*}
\Psi^j\pi(1)&:=&n\\
 \Psi^j\pi(i)&:=&\pi(i-1)\quad  \hbox{for every}\quad 1<i<j\\
 \Psi^j\pi(j)&:=&n-1\\
\Psi^j\pi(i)&:=&\pi(i-2)\quad  \hbox{for every}\quad j<i\leq n.
\end{eqnarray*}
In other words,  the permutation $\Psi^j\pi$ is obtained from $\pi$ by prefixing $n$ in the position preceding the first number in $\pi$ and inserting $n-1$   in the $j$-th position of the resulting permutation.
For a set $A$ of permutations 
we denote by $\Psi^j(A)$ the set of the images by $\Psi^j$ of all the
permutations of $[n-2]$ belonging to $A.$ As a last element of notation, let
us denote by $S^j$ the set of those permutations $\tau$ of $[n]$ for which  
$\tau^{-1}(n-2)<j.$
Consider 
$$A^j:= \Psi^j(B_{n-2})\cap S^j$$
and set 
$$A_n:=\cup_{j=2}^{n}A^j.$$
(The attentive reader may note that $A^2=\emptyset$ but we felt it more natural
not to exclude this set from the above union.)
As every permutation in $A_n$ has $n$ at its first position no two of them
can be cyclic shifts of each other, whence $|B_n|=n|A_n|$.  
Therefore in order to check that we have constructed the right number of
permutations it is sufficient to verify that 
\begin{equation}\label{eq: rn}
|A_n|={{n-1} \over 2}t_{n-2}
\end{equation}
To this effect, recall that by our hypothesis the set $B_{n-2}$ is  invariant
with respect to cyclic shifts.  
This implies that the number of those of its sequences in which a fixed
element, in our case $(n-2),$ is  
confined to any particular subset of the coordinates is proportionate to the
cardinality of the coordinate set in question,  
and thus
$$|A^j|={j-2 \over n-2}|B_{n-2}|,$$
whence
$$|A_n|=\sum_{j=2}^n |A^j |=\sum_{j=2}^n  {j-2 \over n-2} |B_{n-2}|=
{n-1 \over  2}t_{n-2},$$
which, substituting the value of $t_{n-2}$, yields
$$|A_n|={(n-1)! \over 2^{\lfloor{n \over 2}\rfloor}}.$$
This settles our claim (\ref{eq: rn}) and proves that $B_n$ has the requested
number of permutations. 

To conclude the proof it remains to show that every pair of sequences from $B_n$ represents a 
$G(\overline{\{1\}})$--different pair of permutations.
We will first prove that such is the case if both sequences are from $A_n$. 
If they belong to the same 
$A^j$ then this is obvious since the two permutations in such a pair must differ somewhere in those coordinates where they feature 
an element of $B_{n-2}$ and thus the corresponding elements of $B_{n-2}$
must be different sequences. This implies, by our hypothesis, that they differ in some
coordinate by strictly more than 1.  
If the two sequences, $\pi$ and $\tau,$ do not belong to 
the same $A^j$, then we must have, say $\pi\in A^j$ and $\tau \in A^k$ with $j<k.$ But then in the 
$k$-th position $\tau(k)=n-1$, while by definition, $\pi(k)<n-2,$ settling
this case as well. 

If $\pi$ and $\tau$ are two permutations that do not belong to $A_n$ but have
the value $n$ in the same position, then they are clearly in a similar
relation as their respective cyclic shifts in $A_n$, thus the above argument
still applies.

Finally, we must prove that 
any two of our sequences having the symbol $n$ in different positions also
represent a $G(\overline{\{1\}})$--different pair of permutations. Now, unless
the symbol $n$ of both of the two sequences meets the symbol $(n-1)$ of the
other one, we are done. Otherwise they have their respective subsequences
belonging to $B_{n-2}$  
positioned in the very same coordinates and it suffices to see that these
subsequences are different. For this purpose  
suppose that the two sequences have their symbol $(n-1)$ in the $j$-th and the
$k$-th position, respectively.   
 But then, supposing $j<k$ we can say that they must have their respective
 symbols  
$(n-2)$ in different positions since the one having its $(n-1)$ in the $k$-th
position has its $(n-2)$ in a position belonging to the open interval $(j, k)$
while the other one has it in the complement of the closed interval $[j, k]$
by construction. 
This proves our theorem for every odd $n.$ 

In order to prove our claim also for even values of $n,$ it is enough to
consider the set $A_{n+1}$ (now $n+1$ is odd) and delete the first entry,
which is $(n+1)$, from each of the permutations in this set. This way we get  
the right number of permutations of $[n]$ and their pairwise relations satisfy
the requirement by the previous part of the proof. 
\hfill$\Box$

\medskip
\begin{rem}\label{rem:chi}
{\rm Consider the graph whose vertex set
is  the set of permutations of $[n]$  and such that  two permutations form an edge if and only if they satisfy the requirement we dealt with in Theorem~\ref{thm:com}. Denote this graph by
  $H_{\overline{\{1\}}}(n)$. Observe that its clique number
  $\omega(H_{\overline{\{1\}}}(n))=T(n,\overline{\{1\}})$ by definition and
  notice 
  that by the proof above its chromatic number
  $\chi(H_{\overline{\{1\}}}(n))$ has the same value. (The sets $C(\pi)$
  defined in the proof can serve as color classes of an optimal coloring.) 
  This observation will be used in the proof of the subsequent corollary.}
\end{rem}

\medskip

With some additional argument the above theorem gives the exact
value of $T(n,\overline{\{q\}})$ also for $q\neq 1$. 
We will need the following well-known lemma, the roots of which go back to
Shannon \cite{Sh}. We give a short proof for the sake of completeness. 

\begin{lem} \label{lem:perf}
Let $G_1,\dots,G_k$ be graphs and let $G_1\cdot\dots\cdot G_k$ denote their
co-normal product, i.e., the graph with vertex set $V(G_1)\times\dots\times
V(G_k)$ in which  two vertices $\mbf x, \mbf y$ are adjacent if there is an $i$
such that the respective  $i$-th entries  $x_i, y_i$ of these sequences satisfy
$\{x_i,y_i\}\in E(G_i)$. If $\chi(G_i)=\omega(G_i)$ holds for every $i$, then
$\omega(G_1\cdot\dots\cdot G_k)=\prod_{i=1}^k \omega(G_i).$   
\end{lem}

\proof
It is easy to verify that $\omega(G_1\cdot\dots\cdot G_k)\ge\prod_{i=1}^k
\omega(G_i)$ always holds. To prove the reverse inequality observe that
$\chi(G_1\cdot\dots\cdot G_k)\leq \prod_{i=1}^k \chi(G_i)$. By
$\omega(G_1\cdot\dots\cdot G_k)\leq\chi(G_1\cdot\dots\cdot G_k)$ the
conditions $\chi(G_i)=\omega(G_i)$ imply the statement.
\hfill$\Box$

\begin{cor}\label{cor:alt}
Let $q$ be an arbitrary fixed natural number and let $n$ have the form $aq+m$,
where $m\in \{0,\dots,q-1\}$. Then 
$$T(n,\overline{\{q\}})=
\frac{n!}{\left(2^{\lfloor\frac{a}{2}\rfloor}\right)^{q-m}
  \left(2^{\lfloor\frac{a+1}{2}\rfloor}\right)^m}.$$  
\end{cor}

\proof
Let $S_n$ be the set of all permutations of $[n]$ represented as sequences and 
consider a largest possible set $B_n$ of sequences from $S_n$ which satisfies
the requirements for $D=\overline{\{q\}}$.   
Let $h:\N\to \{0,\dots, q-1\}$ be the residue map modulo $q$, or, in fact, any
map for which 
$h(k)=h(\ell)$ if and only if $q$ divides $|k-\ell|$. For sequences $\mbf
x=x_1\dots x_n$ extend $h$ as $h(\mbf x):=h(x_1)\dots h(x_n)$. Partition $S_n$
according to the image of $h$, i.e., put $\mbf x$ and $\mbf y$ into the same
partition class iff $h(\mbf x)=h(\mbf y)$. The number of partition classes so
obtained is 
$$t:=\frac{n!}{(a!)^{q-m}((a+1)!)^m}={n\choose{a,\dots,a,a+1,\dots,a+1}}$$ 
We call the classes $W_1,\dots, W_t$. If two sequences $\mbf x, \mbf y$ belong
to different 
$W_j$'s then there must be a position $i$ for which $|x_i-y_i|$ is not
divisible by $q$, in particular, it is not equal to $q$. Thus
$T(n,\overline{\{q\}})$ is just the sum of the maximum possible cardinalities
of sets of sequences one can find within each $W_j$ such that each pair of
these sequences satisfies the condition. 

Fix any class $W_j$. For each $\mbf x\in W_j$ and each position $i$ the
value $h(x_i)$ is the same by definition. Let $h^j_i$ denote this common
value. For $k\in\{0,\dots,q-1\}$ set $E_k=\{i\;|\; h^j_i=k\}$. Consider the
subsequence of each $\mbf x\in W_j$ given by the entries at the positions
belonging to $E_k$. Note that the size of $|E_k|$ is either $a$ or $a+1$.
Let $H_k$ be the following graph. Its vertex set consists of $|E_k|$-length
sequences of different 
numbers from $[n]\cap \{\ell\;|\; h(\ell)=k\}$. Two such sequences
$\mbf x$ and $\mbf y$ are adjacent in $H_k$ iff at some coordinate $i$ we have
$|x_i-y_i|\neq q$. It is straightforward that $H_k$ is isomorphic to the graph 
$H_{\overline{\{1\}}}(|E_k|)$ defined in Remark~\ref{rem:chi}. Whence its
clique number is
$T(|E_k|,\overline{\{1\}})$, while, by Remark~\ref{rem:chi}, its
chromatic number 
has this same value. Let $\hat H_j$ be the graph with vertex set $W_j$ where
two vertices are adjacent if they satisfy the requirement that at some
position their difference is neither $0$ nor $q$. One easily verifies that
$\hat H_j$ is 
isomorphic to the co-normal product (for the definition see
Lemma~\ref{lem:perf}) of the graphs $H_0,\dots, H_{q-1}$, which 
is, by the foregoing, isomorphic to $\prod_{k=0}^{q-1}
H_{\overline{\{1\}}}(|E_k|)$. We are interested in the clique number of this
graph. By Lemma~\ref{lem:perf} and Remark~\ref{rem:chi} this value is
equal to $\prod_{k=0}^{q-1}\omega(H_{\overline{\{1\}}}(|E_k|))$. Noticing that $q-m$
of the sets $E_k$ have size $a$ and $m$ of them have size $a+1$, 
this is further equal to
$\left({{a!}/{2^{\lfloor\frac{a}{2}\rfloor}}}\right)^{q-m}
\left({{(a+1)!}/{2^{\lfloor\frac{a+1}{2}\rfloor}}}\right)^{m}$ 
by Theorem~\ref{thm:com}.

The latter value is the same for all sets $W_j$ and the number of these sets
is ${n\choose{a,\dots,a,a+1,\dots,a+1}}$ (with $a$ and $a+1$ appearing $q-m$
and $m$ times, respectively). Thus we have obtained
$$T(n,\overline{\{q\}})={n\choose{a,\dots,a,a+1,\dots,a+1}} 
\left({{a!}\over{2^{\lfloor\frac{a}{2}\rfloor}}}\right)^{q-m}
\left({{(a+1)!}\over{2^{\lfloor\frac{a+1}{2}\rfloor}}}\right)^{m}$$
$$= 
\frac{n!}{\left(2^{\lfloor\frac{a}{2}\rfloor}\right)^{q-m}
  \left(2^{\lfloor\frac{a+1}{2}\rfloor}\right)^m}.$$  
\hfill$\Box$

\section{Graph pairs}
It seems interesting to study the relationship of the values of $T(n, D)$ 
for pairs of disjoint sets (graphs) and their union, especially in case of pairs of 
complementary sets. 

Let us define 
$$\phi(D, \overline{D}):=
\limsup_{n\rightarrow \infty}{1 \over n}\log{T(n,D)T(n, \overline{D}) \over n!}$$
and call it the {\it split strength} of the partition $\{D, \overline{D}\}$ of the natural numbers. Consider the case $D:=\{1\}.$ We know from \cite{KMS} that 
$${10}^{n-4\over 4}\leq T(n,\{1\})\leq 2^n.$$
(We do not need the sharper form of the upper bound here. For an exponential improvement in the above lower bound the reader is invited to consult Brightwell and 
Fairthorne \cite{BF}.) 
Using this in combination with Theorem \ref{thm:com} yields

\begin{prop}\label{prop:sp1}
$$0.33 <\phi(\{1\}, \overline{\{1\}})\leq {1 \over 2}$$
\end{prop}

We continue with other examples. Denoting by $2\N$ the set of the even 
numbers, we would like to determine $\phi(2\N, \overline{2\N}).$
To this end, notice first that
$$T(n, \overline{2\N})={n \choose \lfloor {n \over 2}\rfloor}.$$
In fact, this easily follows, as in \cite{KM}, by observing that two
permutations differ in every position by an even number if and only if the
even numbers occupy the same set of positions in both. 

Somewhat surprisingly, $T(n, 2\N)$ seems hard to determine and we only have
some easy  bounds. 

\begin{prop}\label{prop:sp2}
$$\frac {n!\left(\left\lceil\frac{n}{2}\right\rceil +1\right)} {2{n \choose
    \lfloor\frac{n}{2}\rfloor}} 
\leq T(n, 2\N) \leq {n! \over 2^{\lfloor{n \over 2}\rfloor}}$$
\end{prop}

\proof
The upper bound is a trivial consequence of (the upper bound part of) Theorem 1.
Although the lower bound follows from the lower bound on $\kappa(K_n)$ in
\cite{KMS}, yet for the reader's
convenience 
we give the details without explicit reference to said paper. 
(Those needing more details may however consult \cite{KMS}).
We consider the set $[n]$ as the disjoint union of its respective subsets of odd and even numbers. Correspondingly, we divide the coordinate set in two (with a little twist). 
In the first $\lceil \frac{n}{2}\rceil+1$ coordinates we write the even
permutations of 
the set $A$ consisting of all the odd numbers from $[n]$ with the addition of the 
extra symbol $\star.$ (More precisely, first we represent these $\lceil \frac{n}{2}\rceil+1$
many symbols bijectively by the first natural numbers up to their cardinality, then extend this bijection to the permutations of both sets and consider only those permutations of the elements of $A$ that correspond to the even permutations of the first $|A|$ natural numbers).  We represent an arbitrary permutation of $A$ in form of a sequence 
${\mbf x}$ and similarly let $\mbf y$ be an arbitrary permutation of the set $B$ of the even elements 
of $[n].$ We will say that $\mbf y$ is {\it hooked up to} $\mbf x$ if we replace the $\star$ in 
${\mbf x}$ by the first coordinate of $\mbf y$ and concatenate the rest of $\mbf y$ as a suffix 
to the resulting sequence. Let us denote by $\mbf x \leftharpoonup \mbf y$ the  
permutation of $[n]$ so obtained. Define
$A\leftharpoonup B$ to be the set of all these permutations as $\mbf x$ and 
${\mbf y}$ take all of their possible values. Clearly, 
$$|A\leftharpoonup B|=\frac{1}{2}
\left (\left \lceil \frac {n}{2}\right \rceil+1\right )! \left \lfloor
  \frac{n}{2}\right \rfloor!$$ 
which in turn equals the claimed lower bound in the statement of the Proposition. 
It is very easy to see on the other hand that all the pairs of permutations from 
$A\leftharpoonup B$ differ by an even number in some coordinate.
 
\hfill$\Box$
\begin{cor}\label{cor:odd}
$$0\leq\phi(2\N, \overline{2\N})\leq \frac{1}{2}$$
\end{cor}

\hfill$\Box$

Next we quickly review the following immediate consequence of our
hitherto results on split strength.

\begin{prop}\label{prop:d}
Let $q$ be an arbitrary but fixed natural number. Then $\phi(\{q\},
\overline{\{q\}})$ 
is independent of the actual value of $q$.
\end{prop}

\proof
We prove more, namely that the asymptotics of $T(n, \{q\})$ 
is independent of the value of $q$ we fix and the same is true for $T(n,
\overline{\{q\}}).$ For the latter it follows immediately from the
formula given in Corollary \ref{cor:alt}.

Now we turn to $T(n, \{q\}).$ Consider the distance graph $G(\{q\})$ of the
set $\{q\}$ and  
look at the graph it induces on $[n].$ Since the latter is isomorphic to a
subgraph of $P_n,$ the path on $n$ vertices that the analogous distance graph
$G(\{1\})$ induces on the same set, we immediately see that 

\begin{equation}\label{eq:iso}
T(n, \{q\})\leq T(n, \{1\}).
\end{equation}
In the reverse direction, we just have to observe that, for every $m\in
\{0,1,\dots, q-1\}$,  
the graph $G(\{q\})$ induces an infinite path on the residue class 
$q\N+m$ of the numbers congruent to $m$ modulo $q.$ This implies
\begin{equation}\label{eq:mp}
T(n, \{q\})\geq \prod_{m=0}^{q-1}T\left (\left \lfloor\frac{n-m}{q}\right
  \rfloor,  \{1\}\right) 
\end{equation}
by concatenating the respective constructions of permutations for each fixed
$m$. 
Whence it is immediate that $T(n, \{q\})$ and $T(n, \{1\})$ have the same
exponential growth rate. 
\hfill$\Box$

We know very little about split strength and thus there are many questions to
ask.  
Is it always true that $\phi(D, \overline{D})$ is finite and non--negative as
it seems by these examples?  
In order to see the greater picture, we have to look at different kind of
growth rates as well.

\section{Intermediate growth}

So far we have only seen growth rates at an exponential factor away from 
either 1 or $n!\ .$ We intend to show here, however, that 
in between growth rates are also possible. In particular, we will see that 
$T(n, D)$ and $T(n, \overline{D})$ can have essentially the same growth rate,
while their product is still about $n!\ $.

Let $ex(n)$ denote the largest  exponent $s$ for which $2^s$ is a divisor of
$n.$  
We define

\begin{equation}\label{eq:exp}
E:=\{n\;|\; n \in \N, \; ex(n)\equiv 0 (\hbox{mod 2})\}.
\end{equation}

\begin{thm}\label{thm:self}
If $n$ is a power of $4$, then we have
\ 
\begin{description}
\item[(a)]
$$(\sqrt{n})!^{\sqrt{n}}\leq T(n, E)\leq \frac{n!}{(\sqrt{n})!^{\sqrt{n}}},$$
\item[(b)]
$$(\sqrt{n})!^{\sqrt{n}}\leq T(n, \overline{E})\leq
\frac{n!}{(\sqrt{n})!^{\sqrt{n}}}.$$ 
\end{description} 
\end{thm}

\proof
We prove the lower bound part of (a) first. 
It will be convenient to consider the elements of $[n]$ as binary sequences of
length  
$t:=\lceil \log n \rceil$, with each natural number from $[n]$ represented by
its {\it binary expansion}. (Integer parts could be deleted by our assumption
on $n$, moreover, we also know that $t$ is an even number.) In fact, instead of
permuting the $n$ integers in $\{1,\dots,n\}$, now we will permute the $n$
numbers in $\{0,\dots, n-1\}$. With a shift by $1$, the two are obviously
equivalent for our purposes.  
For simplicity, we will index the coordinates of the binary expansions 
from right to left. Hence in particular $m$ is odd if in its binary expansion 
$\mbf x=x_t x_{t-1}\dots x_1$ the rightmost coordinate $x_1$ is 1 and even
else. Let further $\mbf x^{odd}$ and $\mbf x^{even}$ denote the subsequence of
the odd and the even indexed coordinates of $\mbf x$, respectively. Finally,
let $\nu(\mbf x)$ be the smallest (i. e., rightmost) 
index $i$ for which $x_i=1.$ By a slight abuse of notation we will consider
the various subsets of $\{0,\dots, n-1\}$ as subsets of $\{0,1\}^t.$  Quite
clearly, for every $\mbf x\in \{0,1\}^t$ we have $$\nu(\mbf x)=ex(\mbf x)+1$$
and, in particular, $\mbf x \in E$ if and only if
$\nu(\mbf x)\equiv 1$ modulo 2. 
In order to prove the lower bound, let us consider the partition induced  on
$\{0,1\}^t$ (i.e., on $\{0,\dots,n-1\}$) by    
the mapping $f:\{0,1\}^t \rightarrow \{0,1\}^{\frac{t}{2}}$
where 
$$f(\mbf x):=\mbf x^{even} \;\hbox{for every}\; \mbf x\in \{0,1\}^t .$$
(The classes of the partition are the full inverse images corresponding to the
various values of $f.$) 
It follows by construction that 

\begin{equation}
f(\mbf x)=f(\mbf y)\quad \hbox{implies} \quad |\mbf x- \mbf y| \in E
\end{equation}
where by the difference of the vectors $\mbf x$ and $\mbf y$ we mean the
difference in ordinary arithmetics of  
the natural numbers they represent. Indeed, executing the subtraction in the
binary number system we  
are using here one sees that both $\mbf x-\mbf y$ and $\mbf y-\mbf x$ have
their rightmost $1$ in the position where, scanning  
the binary expansions of $\mbf x$ and $\mbf y$ from right to left, we find the
first position in which they differ.  
Now, since $\mbf x^{even}=\mbf y^{even}$ by assumption, the position in
question must have an odd index.  
In other words,  $\nu(|\mbf x-\mbf y|)\equiv 1$ modulo 2. 
For every $\mbf z \in  \{0,1\}^\frac{t}{2}$ we denote by
$S(\mbf z)$ the set of all the permutations of the elements of the full
inverse image $f^{-1}(\mbf z)$ of $\mbf z.$ Thus, by our previous argument,
all these permutations are pairwise $G(E)$--different.
Consider the Cartesian product 
\begin{equation}\label{eq:bic}
C:=\prod_{\msbf z \in \{0,1\}^\frac{t}{2}} S(\mbf z).
\end{equation}
Note that the elements of $C$ are permutations of the numbers in
$\{0,\dots,n-1\}$.  
The above consideration implies that $C$ is a set of pairwise
$G(E)$--different permutations.
Further, observing that for every $\mbf z \in  \{0,1\}^{\frac{t}{2}}$ 
$$|f^{-1}(\mbf z)|=2^{\frac{t}{2}}$$
we have 
\begin{equation}\label{eq:crd}
|C|=(2^{\frac{t}{2}})!^{2^{\frac{t}{2}}}
\end{equation}
proving the lower bound in (a). 
(One might get a somewhat larger set by using the hookup operation instead of
straightforward direct product  
but we do not intend to increase the complexity of the presentation for
this relatively small gain here.)

\medskip

Next we prove the upper bound part of (b).
Notice that the set $C$ we have constructed above has a stronger property than
needed so far. In fact, in every coordinate, the absolute difference of our
permutations is either 0 or else it belongs to $E.$ This observation 
will be the basis for our upper bound. 

Consider the auxiliary graph $H_{\overline E}$ the vertices of which are the
permutations of $[n]$ and two are adjacent if they satisfy the
requirement that at some position they have two numbers such that their
difference is in $\overline E$. Clearly, $T(n,\overline E)=\omega(H_{\overline
  E})$, the clique number of this graph by definition. The above observation
about $C$ implies that its independence number $\alpha(H_{\overline E})$ is at
least $|C|$. Note that $H_{\overline E}$ is vertex transitive (any permuting
of the coordinates in the vertices gives an automorphism), thus by a well-known
fact (cf., e.g., in \cite{SchU}) its fractional chromatic number
$\chi_f(H_{\overline E})$ is equal to the ratio of the number of vertices and
the independence number. Using also that the clique number cannot acceed the
fractional chromatic number (cf. \cite{SchU}) we obtained

$$T(n,{\overline E})=\omega(H_{\overline E})\leq \chi_f(H_{\overline E})=
{|V(H_{\overline E})|\over {\alpha(H_{\overline E})}}\leq {n!\over |C|}=
{n!\over (2^{\frac{t}{2}})!^{2^{\frac{t}{2}}}}$$
proving the upper bound in part (b). 
\medskip

Exchanging the role of even and odd above we obtain the lower bound in (b) and
the upper bound in (a) in a similar way.
\hfill$\Box$
\medskip

Theorem~\ref{thm:self} shows that the investigated values of $T(n, E)$ and
$T(n, \overline{E})$ have the same growth rate at around $\sqrt{n!}\ .$ 
The following statement is a straightforward generalization of the above. 

\begin{thm}
For every rational number $\alpha \in (0, 1)$, there is a set
$E_{\alpha}\subseteq \N$ such that for infinitely many values of $n$ we have 
$${(n^{\alpha})!^{n^{1-\alpha}}}\leq T(n,E_{\alpha})\leq {n!\over
    {(n^{1-\alpha})!^{n^{\alpha}}}}$$ and 
$${(n^{1-\alpha})!^{n^{\alpha}}}\leq T(n,\overline{E_{\alpha}})\leq {n!\over
    {(n^{\alpha})!^{n^{1-\alpha}}}}.$$
\end{thm}

\begin{rem}
{\rm Notice that taking logarithm and using the estimate $\log (k!)\approx
  k\log k$ the above inequalities give that $\log T(n,E_{\alpha})$ is about
  $\alpha\log (n!)$, while $\log T(n,{\overline E_{\alpha}})$ is about
  $(1-\alpha)\log (n!)$.}
\end{rem}

\proof
Let 
$$E_{\alpha}:=\{m\;|\; m \in \N, \; ex(m)\equiv 0,1,\dots,p-1 (\hbox{mod}\; q)\}.$$

Set $\alpha=p/q$ and suppose $n$ is a power of $2^q$. The reasoning is essentially the
same as in Theorem~\ref{thm:self}, which is the case $\alpha=1/2, q=2$. 
Instead of $[n]$ we again permute the
elements of $\{0,\dots,n-1\}$ and represent each of these numbers by their
binary expansion. We collect into one group those numbers of $\{0,\dots,n-1\}$
whose binary expansion has the very same subsequence in those
positions which are indexed by numbers congruent to $p,\dots,q-1$
modulo $q$. There 
are $n^{{q-p}\over q}=n^{1-\alpha}$ different such groups each containing
$n^{\alpha}$ different numbers. Permuting the numbers within a group we
get $(n^{\alpha})!$ permutations of those numbers and these are bound to differ
at some position by the difference of two different numbers in the
group. Such a value belongs to $E_{\alpha}$ by construction. Putting all
permutations of all our groups together we obtain
$(n^{\alpha})!^{n^{1-\alpha}}$ different permutations altogether that not only
satisfy the requirements given by the set $E_{\alpha}$ but no two of which
satisfy the requirements prescribed by $\overline E_{\alpha}$. This gives the
upper bound for $T(n,{\overline E_{\alpha}})$ in a similar way as the upper
bound on $T(n,{\overline E})$ is proven in Theorem~\ref{thm:self}. The rest is
also similar to what we have seen there. 
\hfill$\Box$

\section{Exponential growth and Shannon capacity}
In this section we return to the more familiar territory of distance graphs
with finite chromatic number. The relevance of this parameter is shown by the  
following simple observation.
\begin{prop}\label{prop:chr}
Let $G$ be an infinite graph with finite chromatic number $\chi(G).$ Then
$$T_G(n)\leq (\chi(G))^n$$
\end{prop}

\proof
Let $c:V(G) \rightarrow [\chi(G)]$ be an optimal coloring of the vertices of $G$ and let 
$c_n:V(G)^n \rightarrow [\chi(G)]^n$ be its usual extension to sequences. Notice 
that none  of the full inverse images $c_n^{-1}$ of the elements of $
[\chi(G)]^n$ can contain two pairwise $G$--different permutations of
$[n].$ 

\hfill$\Box$

In particular, distance graphs of "rare" sets of distances have finite
chromatic number. 
More precisely, by a result of Ruzsa, Tuza and Voigt \cite{RTV},  
if the set $D:=\{d_1,d_2,\dots d_n,\dots \}$ has the density of a geometric
progression  
in the sense that $\liminf_{n \rightarrow \infty} \frac{d_{n+1}}{d_n} >1,$
then the distance graph  
$G(D)$ has finite chromatic number. Clearly, this density condition is
sufficient but not necessary for the chromatic number to be finite (cf., e.g.,
the set of odd numbers as differences that result in a bipartite graph.)
 
However, for some graphs $G$ with finite chromatic number one can give a
better upper bound on $T_G(n)$.  
This bound is easily obtained once we realize the
tight connection of our present problem with the classical concept of Shannon
capacity of a graph \cite{Sh}.  

Given a sequence $\mbf x\in V^n$ we shall denote by $P_{\msbf x}$ the
probability 
distribution on the elements of $V$ defined by
	$$P_{\msbf x}(a)={1\over n}|\{i\;|\; x_i=a,\ i=1,2,\ldots,n\}|$$
for every $a \in V$; here $\mbf x=x_1\dots x_n$. The probability distribution
$P_{\msbf x}$ is called the {\it type} of $\mbf x$.
Let $V^n(P,\varepsilon)$ denote the set of those $\mbf x\in V^n$ for which
	$$|P_{\msbf x}-P|=\max_{a\in V}|P_{\msbf x}(a)-P(a)|\leq\varepsilon.$$
Let $G$ be a finite graph with vertex set $V=V(G)$ and edge set $E(G).$
As always, we will say that the sequences $\mbf x=x_1x_2\dots x_n \in V^n$ and
$\mbf y=y_1y_2\dots y_n  \in V^n$ are $G$--different if there is at least one index $i \in [n]$
with  
$\{x_i, y_i\}\in E(G).$
Let $\omega(G,n)$ and 
$\omega(G,P,\varepsilon,n)$ be the largest cardinality of any set $C\subseteq
V^n$ and $C'\subset V^n(P,\varepsilon)$, respectively, of pairwise $G$-different
sequences.  
The Shannon capacity of $G$ (or of $\overline G$ in the more usual notational
convention, 
cf. \cite{Sh}) can be defined as $$C(G)=\limsup_{n\to\infty}\frac{1}{n}{\log
       \omega(G,n)},$$
while the {\it capacity $C(G,P)$ of $G$ within the type $P$} is given (cf. 
\cite{CK}) as
$$C(G,P):=
     \lim_{\varepsilon\to 0}\limsup_{n\to\infty}\frac{1}{n}{\log
       \omega(G,P,\varepsilon,n)}.$$ 
It is immediate from the definitions that $C(G,P)\leq C(G)$ always holds and
using the methods of \cite{CsK} it is easy to prove that in fact $C(G)=\max_P
C(G,P)$. 

In what follows we will restrict attention to graphs we call {\it residue
  graph}s. 
We say that an infinite graph $G$ with vertex set $\N$ is a residue graph if 
there exists a natural number $r$ and a finite graph $M=M(G)$ with vertex set 
$\{0,1,\dots, r-1\}$ such that 
$$\{a,b\}\in E(G)\quad \hbox{if and only if} 
\quad \{(a)_{{\rm mod}\; r}, (b)_{{\rm mod}\; r} \} \in E(M)$$

Let $Q$ be the uniform distribution on $\{0,1,\dots, r-1\}.$ We have

\begin{thm}\label{thm:ty}
$$\lim_{n \rightarrow \infty}\frac{1}{n}\log T_G(n)=C(M(G), Q).$$
\end{thm}

\proof
To prove 
$$\limsup_{n \rightarrow \infty}\frac{1}{n}\log T_G(n)\ge C(M(G), Q)$$
consider, for every $n$ those 
sequences  $\mbf x \in \{0,1,\dots, r-1\}^n$ whose type $Q_n$ satisfies
$$Q_n(a)=\frac{1}{n}|\{m \;| \;m\leq n,  (m)_{{\rm mod}\; r}=a\}|$$
for every $a\in \{0,1,\dots, r-1\}.$
Let $M^n$ be the graph whose vertices are the $n$--length sequences 
of vertices of $M$ and whose vertices are adjacent if the corresponding
sequences are $M$--different. 
Let $C_n$ be a complete subgraph of maximum cardinality $M^n$ induces on 
the set of sequences of type $Q_n.$
To any sequence $\mbf x \in V(C_n)$ we associate a permutation of $[n]$ by 
replacing the occurrences of $a$ in the sequence by the different numbers 
congruent to $a$ modulo $r$, in a stricly increasing order. The result is a set
of permutations which is $G$-different. 

For the reverse inequality let $n$ be a multiple of $r$ and consider any
construction achieving $T_G(n)$, i.e., a set of 
permutations of the elements of $[n]$ that are pairwise $G$-different, while
the cardinality of the set is $T_G(n)$. Substitute the occurence of the
number $i$ in each of these permutations by the unique $j\in\{0,\dots,r-1\}$
which is congruent to it modulo $r$. Doing this for all $i\in [n]$ we get
$T_G(n)$ different sequences of vertices of $M$ each having type $Q$ that actually 
form a clique in $M^n$. 
\hfill$\Box$
\medskip

The following corollary is immediate. 

\begin{cor}
$$\lim_{n \rightarrow \infty}\frac{1}{n}\log T_G(n)\leq C(M(G)).$$
\hfill$\Box$
\end{cor}

\medskip\noindent
Let us consider the following 

\par\noindent{\bf Example}
Let $G$ have vertex set $\N$ and set
$$\{a, b\}\in E(G) \quad \hbox{if}\quad |a-b|\equiv 1\;{\rm or}\;
4\;(\rm{mod}\; 5).$$ 
As an easy consequence of Lov\'asz' celebrated formula \cite{LL} for the 
Shannon capacity of the pentagon graph we obtain, using the last theorem (and
also that the Shannon capacity of $C_5$ is obtained by sequences the type of
which is the uniform distribution), that 
$$\lim_{n \rightarrow \infty}\frac{1}{n}\log T_G(n)=\frac{1}{2}\log 5.$$

\bigskip

It is an easy observation that for any graph $M$ and any rational probability
distribution $P$ on its vertex set one can construct (by simply substituting
each vertex by an independent set of appropriate size) a graph $\hat M$ for
which $C(M,P)=C(\hat M,Q)$, where $Q$ is again the uniform distribution. It is
then easy to construct an infinite graph $G$, which is a residue graph with
respect to $\hat M$ and thus the asymptotics of $T_G(n)$ is in an analogous
relationship with the capacity $C(M, P)$ as the one expressed in
Theorem~\ref{thm:ty}. Taking into account the remark that $C(M)$ can be
expressed as the maximum of the values $C(M,P)$ over $P$, we can conclude that
the class of problems asking for the asymptotics of $T_G(n)$ for various
infinite graphs $G$ contains the Shannon capacity problem of all such graphs
for which Shannon capacity is attained as the capacity within a type for some
rational distribution.

\end{document}